# MultiResolution Anomaly Detection Method for Long Range Dependent Time Series

## Lingsong Zhang

*Department of Biostatistics*
*Harvard School of Public Health*
*655 Huntington Ave*
*Boston, MA, 02115.*
*e-mail:* zhang@hsph.harvard.edu

## Zhengyuan Zhu and J. S. Marron

*Department of Statistics and Operations Research*
*University of North Carolina at Chapel Hill*
*Chapel Hill, NC, 27599-3260*
*e-mail:* zhuz@email.unc.edu, marron@email.unc.edu

**Abstract:** Driven by network intrusion detection, we propose a MultiResolution Anomaly Detection (MRAD) method, which effectively utilizes the multiscale properties of Internet features and network anomalies. In this paper, several theoretical properties of the MRAD method are explored. A major new result is the mathematical formulation of the notion that a two-scaled MRAD method has larger power than the average power of the detection method based on the given two scales. Test threshold is also developed. Comparisons between MRAD method and other classical outlier detectors in time series are reported as well.

## 1. Introduction

Detection of Internet intrusions, for example, a Distributed Denial of Service (DDoS) attack, is an important problem in network research and engineering. In a DDoS attack, an attacker may use the security weaknesses (of computers or network protocols) to control other people's computers. The attacker then uses those distributed computers as "zombie" armies to send a large amount of data to one particular website or a particular sub-net, so that legitimate users can no longer access or use the network resources (CERT-Coordination-Center, 2001). Various detection methods have been developed in the field of intrusion detection (see for example Lee et al. (1999); Paxson (1999); McHugh (2001); Barford et al. (2002); Lee et al. (2003)). One important class of methods views the intrusions as a type of network anomaly (McHugh, 2001). This type of method assumes that, the typical traffic, such as browsing websites, downloading files, etc, is the "normal" network traffic. Malicious usages including network intrusions are treated as network anomalies. This separation of the observed data is similar to the separation of outliers from regular observations (Hawkins, 1980; Barnett and Lewis, 1994). Thus, developing and utilizing statistical outlier detection methods may help to identify network anomalies.

Network features, such as packet count, byte count and flow count, collected at a single location, form a time series with the following two special features. The first one is Long Range Dependence (LRD) (Leland et al., 1994), which means that the AutoCovariance Function (ACF) $\gamma(k)$ decays at a polynomial rate, i.e., much slower than the usual exponential decay, such as in the classical AutoRegressive Moving Average (ARMA) time series. The second feature is Self-Similarity (SS) (Willinger et al., 1997), i.e., the ACF of the time series at different scales are essentially the same (see formal definitions of LRD and SS in Section 3.1). Barford et al. (2002) showed that different types of network anomalies lead to statistically abnormal signals at different time scales. The LRD and SS properties for the "normal" traces, and the multiscale property of the network anomalies, motivate us to develop a novel detection method, which efficiently uses these multiscale properties.

In the time series outlier detection context, Fox (1972) introduced the notions of additive outlier and innovation outlier. Classical outlier detection methods include intervention analysis (Tiao and Tsay, 1983; Chang et al., 1988; Tsay, 1988), robust methods (Martin and Yohai, 1986) and dynamic models (West et al., 1985). The intervention model was first introduced by Box and Tiao (1975). These classical methods usually assume the observations are independent or short range dependent (i.e., the ACF decays exponentially as the lag goes to infinity). For detecting outliers in LRD time series, these methods may not be suitable, because the artifacts that are naturally generated by the LRD may cause these methods to misidentify more regular observations as outliers (i.e., increased false alarms). Note that identifying some type of outliers, such as the level shift (defined in Section 2.1), can also be viewed as finding change points in a (stationary or non-stationary) stochastic process or (nonparametric) regression curve. See e.g. Pollak and Siegmund (1985), Müller (1992), and Carlstein et al. (1994), etc.

In Zhang et al. (2008), we proposed a MultiResolution Anomaly Detection (MRAD) method, which uses simple aggregation methods to form multiple-scale time series. The new method identifies anomalies at one specific time location, by using all the observations across the time scales. This method is motivated by the scale space ideas in the curve estimation literature (Chaudhuri and Marron, 2000). A brief introduction of this method is in Section 2, along with a motivating example. In this paper, we focus on the theoretical properties of this MRAD method. When the normal network trace is long range dependent, and the outlier has the format of a local mean level shift, we prove that the power of the MRAD method is larger than the average power at single scales (Theorem 1 in Section 3.2).

The MRAD method we propose is similar to other wavelet-based methods,

such as change point detection methods (e.g. Wang (1995, 1999)), time series outlier detection methods (e.g. Bilen and Huzurbazar (2002)), and other related methods such as Thottan and Ji (2003) etc. The major difference between our work and these similar works includes the follows: (1) Our method can be viewed as a non-standard wavelet method. Unlike the standard wavelet methods, we utilize the long range dependent structure to define the aggregation, and our test is based on the approximation components of a (non-standard) wavelet transform instead of the wavelet coefficients. (2) Our work is motivated by scale space ideas, which view the same phenomena at different scales. These different scales usually are correlated. However, the usual wavelet-based method decomposes the original data into several components, where each component corresponds to a different scale, and different scales usually are uncorrelated. (3) We consider the multiple comparisons across scale at one particular time location. Although other methods use multiple scales as well, they usually consider the multiple comparisons across time at one particular scale. In this paper, we will focus on the statistical outlier method in the time series framework. We will not compare our work with the methods of detecting change points, and leave that as potential future work.

The remaining part of this paper is arranged as follows. A motivating example is shown in Section 2, along with a formal definition of the test problem, aggregation methods, and the test procedure. A visual device for reporting the results, the MRAD outlier map, is also introduced in this section. Section 3 discusses several important theoretical properties of the MRAD method, including the power comparisons between the MRAD method and other single-scale methods. Section 4 discusses the test threshold of this method. Section 5 compares the MRAD method with the classical intervention analysis. Discussion and future work is reported in Section 6.

## 2. A motivating example and related visualizations

In this section, we use a real network trace to illustrate the MRAD method and the related outlier map. The trace is a byte-count time series, which was collected every milli-second (*ms*) from 9:30 pm, Monday April 08, 2002, at the main Internet link between UNC and outside. The estimated Hurst parameter is 0.877, which indicates strong LRD. Figure 1 displays a part (6 minutes) of the entire time series. Besides some of the obvious spikes within the series (e.g. the spike around 300,000), it is hard to find other (possible) anomalies.

The MRAD method includes the following steps:

1. Form multiple-scale time series;
2. For each observation at each scale, determine whether it is an anomaly or not;
3. Report, visualize and interpret the test results.

In this section, we will define the test problem (subsection 2.1), the aggregation method used in this paper (subsection 2.2), the test method (subsection 2.3), the related MRAD outlier map (subsection 2.4) using this motivating example.

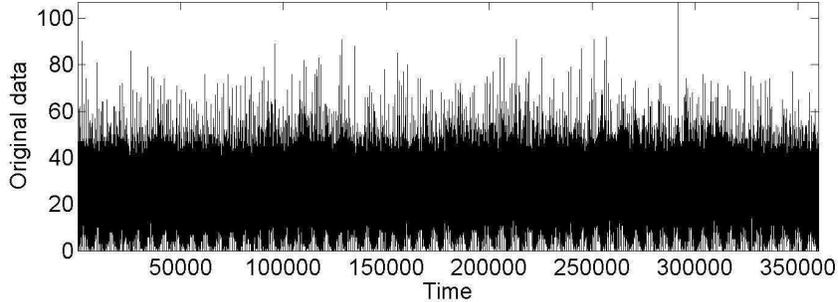

FIG 1. *The byte-count time series, was collected at 9:30 pm, Monday April 08, 2002. The estimated Hurst parameter of the time series is 0.887, which shows long range dependence.*

### 2.1. *The testing problem*

The test problem we intend to solve in this paper can be formulated as the following. Let $\{Y_1(i)\}$, $i = 1, \cdots, N$ be the (standardized) observed time series, e.g. the byte-count time series showed in Figure 1. Because the time series is long range dependent, here we assume the underlying model for $\{Y_1(i)\}$ is

$$Y_1(i) = X_1(i) + \delta I_{i \in [a_0, a_1]}(i), \tag{1}$$

where $\{X_1(i)\}$ is a fractional Gaussian noise (fGn) with Hurst parameter $H$ (see Section 3.1 for formal definitions of fGn and the Hurst parameter). This assumption is made partly for tractability of the theoretical properties (see section 3.2). The other reason that we use fGn as the background trace is that it is a good approximation to highly aggregated Internet traffic time series (e.g. Paxson (1997)).

In the context of intrusion detection, $X_1(i)$ can be viewed as the normal traffic, and the level shift represents a type of network anomaly, such as a DDoS attack. We are interested in detecting the starting time of the attack $a_0$, which can be formulated as a testing problem, i.e., testing whether the observation $Y_1(i)$ is an outlier or not:

$$H_0 : \mathcal{L}(Y_1(i)) = \mathcal{L}(X_1(i)) \text{ vs. } H_1 : \mathcal{L}(Y_1(i)) = \mathcal{L}(X_1(i) + \delta), \tag{2}$$

where $\mathcal{L}(Y_1(i))$ means the distribution of the random variable, $Y_1(i)$. Note that the problem (2) is a pointwise testing problem in the time space.

### 2.2. *Simple aggregation*

There are many methods to form multiple-scale time series, including wavelet methods and kernel methods. We discuss two simple aggregation methods used in this paper, because they are natural, have good theoretical properties (Section 3.2), and are rather easy to be adapted into realtime detection algorithms (See details in Section 3.3.2 in Zhang (2007)). For easy presentation, the multi-scale time series will be constructed by using a dyadic-like structure, which is motivated by the Haar Wavelet bases (see e.g. Ogden (1997) and Vidakovic (1999) for introduction of Haar Wavelet). In this way, the window sizes at different scales increase exponentially. For example, the size of aggregation window at scale $k$ is $2^{k-1}$, i.e., the observation at scale $k$ is a function of consecutive $2^{k-1}$ observations at scale 1. In general, we can replace the base 2 in this dyadic structure to be an arbitrary base $b$, which yields a more general aggregation method.

Assume the same model (1), and the background fGn has Hurst parameter $H$. Let $Y_k(i)$, $i = 1, \cdots, \lceil N/2^{k-1} \rceil$ be the corresponding $k$-scale time series, where $\lceil x \rceil$ returns the smallest integer which is greater than or equals to $x$. One type of aggregation is generated by dividing the original time series into non-overlapping sections, and then aggregate the observations within the sections. Formally, the scale $k$ time series is defined by

$$Y_k(i) = \sum_{j=1}^{L_k} \frac{1}{(L_k)^H} Y_1((i-1)2^k + j),$$

where $L_k = 2^{k-1}$. These special weights are chosen to make $\{Y_k(i)\}$, for all $i$ and $k$, share the same marginal distribution, when there are no outliers in $\{Y_1(i)\}$ (see Section 3.2 for details). By the above definition, the observation $Y_1(i)$ on the finest scale, can only be used once to form the time series at scale $k$. The observations at scale $k$ are a function of the observations within each section. We call this type of aggregation as *Non-Overlapping Window Aggregation* (NOWA). This type of aggregation can be viewed as the approximation component of a nonstandardized wavelet. So our method can be extended to use other (discrete) wavelet methods.

Another type of aggregation is to use overlapping window aggregation. The observation at time $i$ and scale $k$ is defined as

$$Y_k(i) = \sum_{j=0}^{L_k-1} \frac{1}{(L_k)^H} Y_1(i-j),$$

where $i = L_k, L_k + 1, \cdots, N$. This method defines a window of size $L_k$, and then slide it along the original time series, one observation each time, to form all observations at scale $k$. Thus we called it as the *Sliding Window Aggregation* (SWA).

Overlapping window aggregations are widely used, for example, for traditional kernel methods (e.g. Wand and Jones (1995)). Note that the SWA method defined above can be viewed as a special case of the one-sided kernel method (see e.g. Gijbels et al. (1999)). Thus, SWA can be generalized to allow use of a general one-sided kernel (as apposed to the uniform kernel used here). Thus, the usual kernel method can be adapted in this framework for outlier detection.

To identify outliers, classical outlier-detection methods in time series often deal with the whole time series, to estimate the model and then perform hypothesis tests on the residuals. Thus, for example, for the $i$th observation, methods such as the intervention models exploit future information after time $i$. The NOWA method above is of this type. However, detection methods based on realtime cannot use future observations as inputs. Note that the SWA method uses only up-to-now observations, which can be easily adapted for use in real-time detection (see Section III(C) in Zhang et al. (2008)).

These two aggregation methods have a strong relationship. For example, at the $2^k$th time slot, the observation vector $(Y_1(2^k), Y_2(2^k), \cdots, Y_k(2^k))$ are the same for both methods. Thus, the test at each location designed for NOWA is very similar to that for SWA. In fact, the theoretical properties we prove in Zhang (2007) are mainly based on SWA, while all the visualizations in this paper are based on NOWA, because the outlier map based on NOWA method is easy to interpret (see in Section 2.4).

### *2.3. MRAD procedure*

Take one observation as an example, our target problem is to test whether the (standardized) $Y_1(i)$ is an outlier or not. One simple method is to set up the rejection region as $|Y_1(i)| > C_\alpha$, at a given significance level $\alpha$. Our MRAD method uses $\{Y_k(i)\}$, $k = 1, \cdots, M$, to flag outliers. If any of these observations is larger than the threshold (which depends on the aggregation level $M$), we will flag $Y_1(i)$ as an outlier. Formally, the procedure uses $\max_k |Y_k(i)|$ as the test statistic, and $\max_k |Y_k(i)| > C_\alpha^M$ as the rejection region, i.e.,

1. Set a unique threshold $C_\alpha^M$, so that, under the null hypothesis (i.e., there is no outlier in the time series), $P(\bigcup_k \{|Y_k(i)| > C_\alpha^M | \delta = 0\}) = \alpha$.

2. For any $i = 1, \cdots, N$, if any of the observations at some scales exceeds the test threshold, i.e., $\max_k |Y_k(i)| > C_\alpha^M$, we claim that $Y_1(i)$ is a possible outlier.

It can be shown that $C_\alpha^M$ does not rely on $i$, and any scale index $k$, but does depend on the total number of scales. The test threshold will be discussed more intensively in Section 4.

The MRAD method proposed in this particular paper assume that the Hurst parameter is pre-known or can be estimated correctly even when the data trace contains anomalies. There are a lot of literature about the estimation of the Hurst parameter, see for example in Stoev et al. (2005) and references therein. Shen et al. (2007) also provides a robust method to estimate the Hurst parameter in presence of outliers. We will not discuss estimating the Hurst parameter in this paper. We suggest to use a training trace to estimate the Hurst parameter, and then assume the correlation structure of the time series preserves in the near future. Then we can use the estimate and the MRAD method to identify network anomalies. For real-time anomaly detection, we will update the estimate

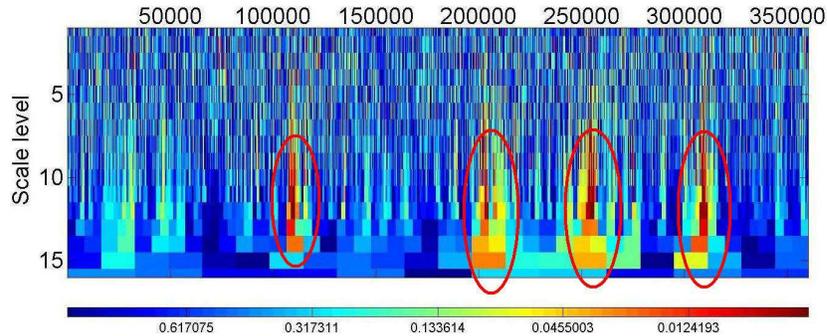

FIG 2. *The color MRAD outlier map for the byte-count time series. Several interesting regions of hot colors are marked by circles.*

of the Hurst parameter periodically (e.g. every day) to improve the detection performance.

### 2.4. MRAD outlier map

In this subsection, we show several visualization devices for the MRAD method - the MRAD outlier map, and the related sliding movie. If there is no anomaly in the original time series, the marginal distributions of all the observations in the multiple-scale time series are the same as $N(0,1)$. Thus, after getting the multiple-scale time series, we can visualize the probabilities of all the observations across time and scales to be outliers, which is the MRAD outlier map.

Figure 2 shows the color MRAD outlier map based on the NOWA method. In the map, the rows report the test results of the time series at different scales, and the columns are the time stamps. Here the scale increases exponentially, i.e., one observation at the $k$th scale contains $2^{(k-1)}$ observations at the finest scale (scale 1). Each column corresponds to the same time location across scales. Each cell of the map display the significance probability ($p$ value) of each observation at a given scale and time location. Within this map, hotter colors (red) represents smaller $p$ values, i.e., the observation has higher chance to be an anomaly; and, cooler colors (blue) corresponds to larger $p$ values, i.e, they are less likely to be anomalies.

From the map, we can find several interesting regions: scale 9-13 at time 100,000; scale 7-15 at time 200,000; scale 7-12 at time 250,000; and scale 7-15 at time 310,000. Note that these locations are not obvious in the original time series (see Figure 1). These regions might correspond to some hidden network anomalies.

There will be an overplotting problem (i.e., the resolution for this image is not enough to visualize the results of all locations) for the above MRAD outlier map. It is natural to zoom the above map and show the test result locally. We develop a sliding movie to display the time series and the outlier map locally and dynamically. Figure 3 shows one carefully selected snapshot of the MRAD outlier movie. It is a zoomed version of the map in Figure 2, where the time range is from 33,000 to 36,000. The map shows two outlying regions, one is around time location 35,600 and scales 5-8; and the other is around time location 33000 and scales below 8. The first one corresponds to a transitional mean level shift, and the second one was caused by no connection for a short period of time. More snapshots of the movies, and related matlab functions can be accessed at the website of Zhang (2006).

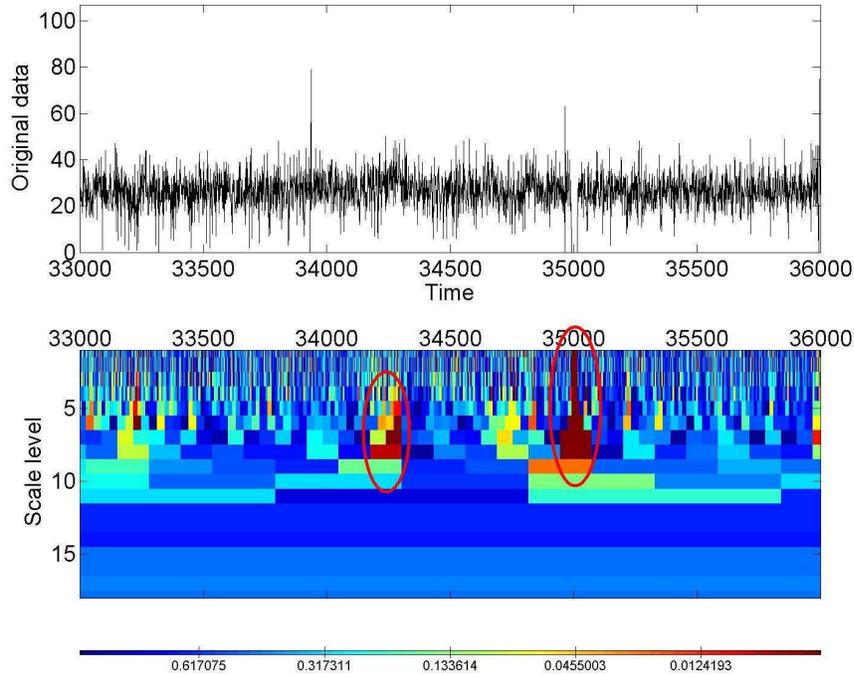

FIG 3. *One carefully selected snapshot of the sliding MRAD outlier movie for the motivation example. The left red block corresponds to a local mean level shift, and the right flagged anomalies were caused by no connection at that specific range of time.*

The outlier map is motivated by the SiZer map in the curve estimation context (Chaudhuri and Marron, 1999). It is also similar to the scalogram in the wavelet context (see in Mallat (1999)). All these graphs are multiresolution visualizations. Note that the major difference between our outlier map and other related visualizations is that we plot the *p*-values of some tests, which can be viewed as a function of the approximation component of a (non-standard) wavelet transform. The SiZer map plots the decisions of some particular tests at a pre-defined significance level (e.g. $\alpha = 0.05$). The scalogram is a graph of the energy of the coefficients of the wavelet transform.

### *2.5. Some remarks on the Hurst parameter*

## 3. Theoretical properties for the MRAD method

In this Section, we study theoretical properties of the MRAD method. Section 3.1 provides a short introduction of fGn and LRD. Section 3.2 discusses some theoretical properties for a two-scale MRAD procedure. We conjecture these properties for two-scale MRAD procedure hold for general MRAD procedures as well.

### *3.1. Background on LRD and fractional Gaussian noise*

Let $\gamma(h) = EX_iX_{i+h}$ be the ACF of a statistical time series. A stationary time series is said to have *long range dependence* (LRD), if $\gamma(h) \propto L(h)h^{-\alpha}$ as $h \to \infty$, where $L(h)$ is a slowly varying function, and $\alpha \in (0,1)$ (Taqqu, 2003).

An important feature of LRD time series is that $\sum\limits_{h=1}^{N} |\gamma(h)| \to \infty$ as $N \to \infty$.

Here we give the definition of fractional Brownian motion (fBm) and fractional Gaussian noise (fGn).

**Definition 1.** *A stochastic process $\{B_H(t)\}_{t \in \mathcal{R}}$ is called a fractional Brownian motion (fBm), if it is a Gaussian process with mean 0, stationary increments, variance $EB_H^2(t) = t^{2H}\sigma^2$ and covariance $EB_H(s)B_H(t) = (s^{2H} + t^{2H} - |s - t|^{2H})\sigma^2/2$, where $0 < H < 1$ is called scaling exponent or Hurst parameter.*

**Definition 2.** *The increment process of a fractional Brownian motion, $X_i = B_H(i+1) - B_H(i)$, $i \geq 1$ is called a fractional Gaussian noise (fGn).*

Note that fGn is a mean zero, stationary Gaussian time series, with AutoCovariance function $\gamma(h)$ given by $\gamma(h) = \{|h+1|^{2H} - 2|h|^{2H} + |h-1|^{2H}\}\sigma^2/2$, $h \geq 0$. For $H \neq 1/2$, $\gamma(h) \sim \sigma^2H(2H-1)|h|^{2H-2}$ as $h \to \infty$. So when $1/2 < H < 1$, the fGn shows long range dependence. It has been used for modeling network traffic (Paxson, 1997). See Mandelbrot and Van Ness (1968); Taqqu (2003) for more information about fractional Brownian motion and fractional Gaussian noise.

### *3.2. Theoretical properties of MRAD procedure*

In this section, we give theoretical results for MRAD. In particular, we show that the power of a two-scale MRAD procedure (either based on NOWA or SWA), as described in Subsection 2.3, is larger than the average power of the naive outlier detection method based on one of the two scales. This result gives theoretical justification for using the MRAD method over methods based on a single scale. Proofs of the propositions and theorems are included in the appendix.

Since $\{X_1(i)\}_{i=1,\cdots,N}$ is a fGn with Hurst parameter $H$, it can be shown that $\{X_L^N(i)\}$ defined by the NOWA aggregation, $X_L^N(i) = \sum_{j=1}^{L} X_1((\lceil i/L \rceil - 1)L + j)/L^H$, is also a fGn with the same Hurst parameter $H$. Thus, $X_L^N(i)$ has the same marginal distribution as $X_1(i)$. The $\{X_L^S(i)\}$ based on SWA, has the same marginal distribution, as the $\{X_L^N(i)\}$ based on NOWA. So for both aggregation methods, all $\{X_L(i)\}$ share the same marginal distribution (i.e. $N(0,1)$). We are interested in testing whether the $i$th observation $Y_1(i)$ is an outlier (i.e. $Y_1(i) = X_1(i) + \delta$, and $\delta \neq 0$) or not, as defined in Equation (2).

Let $\{|Y_L(i)| > C_{\alpha,L}\}$ be the rejection region for scale $L$ for the significance level $\alpha$. Because the marginal distributions of $\{Y_L(i)\}$ are the same for different scales, when there is no outlier, the threshold $C_{\alpha,L}$ does not depend on $L$. We denote this threshold as $C_\alpha$. Note that for a given significance level $\alpha$, $C_\alpha = \Phi^{-1}(1 - \alpha/2)$.

The following propositions and theorems show some important theoretical properties for the MRAD procedure. Unless otherwise specified, all the following results are derived based on SWA (similar results also holds for NOWA).

**Proposition 1.** *If $P_0(\cup_L \{|Y_L(i)| > C_\alpha^M\}) = \alpha$, and $P_0(|Y_L(i)| > C_\alpha) = \alpha$, we have $C_\alpha^M \geq C_\alpha$.*

*Remark*: This proposition shows the MRAD method provides a more conservative threshold than the naive outlier detection method at any scales.

**Proposition 2.** *The time series at the $i$th location over different scales, $\{Y_k(i)\}$, $k = 1, 2, \cdots$, based on the SWA with base $b$, is a stationary process, with ACF*

$$\rho_b(k-1) = \frac{1}{b^{kH}} \left[ 1 + \frac{1}{2}(b^{k2H} - (b^k - 1)^{2H} - 1) \right]$$

$$= Hb^{k(H-1)} + \frac{b^{-kH}}{2} - \frac{H(2H-1)}{2} b^{k(H-2)} + o(b^{k(H-2)}).$$

*Remark*: This proposition shows that at each time location, the observations across scales form a stationary process. The AutoCovariance function decays exponentially, i.e., this process over scales is short range dependent. This also shows that when $H$ is large, the decay rate will be slow.

For the remaining part of this section, let $C_\alpha^M$ be the 2-scale MRAD testing threshold, and $C_\alpha$ be the testing threshold based on one single scale.

**Proposition 3.** *For a 2-scale MRAD method, let $C_\alpha^M$ be the testing threshold of significance level $\alpha$, we have*

$$C_\alpha^M = C_0 - \frac{\phi(C_0)C_0^2 H^2}{2\sqrt{1-\alpha}} L^{2(H-1)} + o(L^{2(H-1)}),$$

*as $L \to \infty$. Here $C_0 = \Phi^{-1}((1 + \sqrt{1-\alpha})/2)$.*

*Remark*: The proposition implies that $C_0$ is the limit of $C_\alpha^M$. When $L$ is large, the testing threshold $C_\alpha^M$ is close to $C_0$. In addition, the convergence rate is of the order $2(H-1)$. Larger $H$ corresponds to lower convergence rates.

**Theorem 1.** *Let* $\beta_{(1,L)} = P_1(\max_{l=1,L} |Y_l(i)| > C_\alpha^M)$, $\beta_1 = P_1(|Y_1(i)| > C_\alpha)$, *and* $\beta_L = P_1(|Y_L(i)| > C_\alpha)$. *For any* $\delta > 0$, *there exists* $\alpha_\delta > 0$ *and* $L_\delta > 0$, *when* $\alpha \in (0, \alpha_\delta)$ *and* $L > L_\delta$, *the following inequality holds:*

$$\beta_{(1,L)} \geq \frac{\beta_1 + \beta_L}{2}.$$

*Remark*: This theorem shows that for any level shift, when the significance level $\alpha$ is small and $L$ is large, the power of the two-scale MRAD is larger than the average power of the outlier detection methods each based on a single scale. In the context of network anomaly detection, the network anomalies in general can exist at any scale (Barford et al., 2002). To detect one particular anomaly, if by chance we use an appropriate scale, the one scale method will have the maximum detection power. In practice, it is hard to know which scale to use. This theorem suggests that the multiresolution method on average will have larger power than the methods based on single scale, which shows the usefulness of the multiscale ideas.

## 4. Test thresholds

We used the asymptotic independence between the maximum value and the minimum value of the above stationary process (as described in Proposition 2), and developed an asymptotic test threshold for an $m$-scale MRAD method, $C_I^M = \Phi^{-1}((1-\alpha)^{1/2m})$. See details in Zhang (2007). Note that this threshold is not precise, because it is the threshold when $m$, the number of scales, goes to infinity. Often a rather small number of scales (e.g. 15) is used, thus, the asymptotic test threshold might be too conservative, and have relatively low power.

We also developed an improved test threshold using computer simulation. The method also uses the result of Proposition 2: at one particular time, the time series across scales is stationary, and the ACF of these stationary processes, at different time locations, are the same. Thus, we only need to simulate once to get the threshold for all locations. We developed a MATLAB function, `mradtestthreshold.m` (available at Zhang (2006)), to approximate the exact threshold based on simulation of multivariate Normal (given the Hurst parameter, the significance level, and the number of scales). We refer to this threshold as *the improved test threshold* for the MRAD method.

Zhang (2007) provided a list of test thresholds at different combinations of the significance level, Hurst parameter, and the number of scales. We do not report them in this paper to save space. We find that at a given number of scales, and a given significance level, the test threshold decreases as the Hurst parameter increases. On the other hand, if the Hurst parameter and the significance level are fixed, the threshold increases when the number of scales increases. In fact, it will approach to the asymptotic threshold as the number of scales goes to infinity. See more discussion in Zhang (2007).

## 5. Comparison between the MRAD method and other related methods

As mentioned in Section 1, there are several different methods for detection of outliers in a time series. In this section, we compare our MRAD method with the classical intervention analysis method. Section 5.1 describes the evaluation metrics we use to compare our method with the intervention analysis methods. Section 5.2 provides the comparisons among them. Section 5.3 provides a semi-experiment to illustrate the usefulness of the MRAD method.

### *5.1. Metrics for evaluation*

The following metrics are used to evaluate and compare our MRAD method with the classical intervention analysis method: the False Discovery Rate (FDR), False Negative Rate (FNR) and the True Discovery Rate (TDR). Here the term "discovery" means declaring an observation as an outlier (positive), and "negative" means declaring an observation as a regular observation. Table 1 is the classical FDR definition table, which was introduced in Benjamini and Hochberg (1995). Here, $U$ is the number of regular observations flagged correctly as "normal traffic", i.e., true negatives; $V$ is the number of regular observations wrongly classified as anomalies, i.e., false positives; $T$ is the number of outliers wrongly declared as "normal traffic", i.e., false negatives; $S$ is the number of outliers correctly flagged as anomalies, i.e., true positives, and $R$ is the total number of observations identified as outliers.

Table 1
*The classical FDR definition table.*

|  | Declared non-outlier | Declared outlier | Total |
|---|---|---|---|
| regular observations | $U$ | $V$ | $m_0$ |
| true outliers | $T$ | $S$ | $m - m_0$ |
|  | $m - R$ | $R$ | $m$ |

Base on the above notations, the above three metrics can be defined formally as follows: *TDR*, $E(S/(m - m_0))$, is the average proportion of declared outliers among the true outliers. *FDR*, $E(V/R)$, is the average false discovery rate, i.e., among all the declared outliers, the average ratio of those that were regular observations. *FNR*, $E(T/(m-R))$, is the average false negative rate, i.e., among all the observations declared not outliers, those are true outliers.

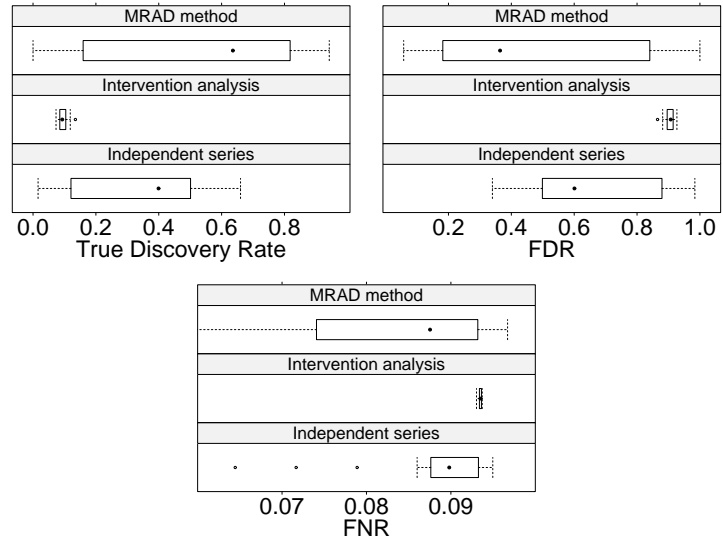

FIG 4. *True discovery rate (the left panel), FDR (the middle panel) and FNR (the right panel) comparisons among the MRAD method (the top row in three panels), the naive one-scale method (the bottom row in three panels), and the intervention analysis method (the middle row in three panels). The MRAD method has the largest median TDR (smallest FDR and FNR), and the median TDR for the intervention analysis is the smallest (largest FDR and FNR). This suggests that our method is better than the intervention analysis, when the time series is long range dependent. Note that the intervention model actually specifies a wrong model.*

### 5.2. *Comparisons between the MRAD method and the intervention analysis*

In this subsection, we use simulation to compare the proposed MRAD method and the classical outlier detection method based on intervention analysis.

The intervention analysis method was first introduced by Box and Tiao (1975). Chang et al. (1988); Tsay (1988) etc. developed an iterative algorithm for detecting outliers in time series based on intervention analysis. The SAS/ETS package (in SAS 9.13) provides routines for detecting outliers in time series, by using the intervention analysis methods. In the following comparisons, we use the SAS/ETS package to detect the imputed level shift automatically.

In this subsection, we simulate the fractional Gaussian noises with $2^{15}$ observations (which roughly approximates a one-hour trace recorded at $100ms$ time intervals). The level shifts we imputed here are simulated from the following distributions: the starting point is from $U[0, 2^{14}]$, the duration of the level shift is from exponential distribution with mean 4000, and the mean level shift is 1 for all simulation, i.e., the intensity is the same as the standard deviation of the background. Note that the uniform distribution of the starting point, and the exponential distribution of the duration, are just the natural starting point to explore and compare the MRAD method and other methods. Exploration of the distributions of the duration and starting time of a network intrusion is potential future work.

In the following plots, we report the above three metrics, FDR, FNR, and TDR, for three types of method: our MRAD method, a naive one-scale method, and the intervention analysis method. The left panel in Figure 4 shows the distribution (by using boxplot) of the True Discovery Rate for these three detection methods. The top panel provides the boxplot for the TDRs using the MRAD method based on 10 sets of simulation. Each simulation set contains 100 simulation of background traces. We calculate the true discovery proportions for each simulation, and show the average of these proportions. The boxplot visualizes the distribution of 10 averages of these proportions. Note that the median TDR of the MRAD method is larger than 0.6. However, the spread of the TDRs is large. The middle panel is the boxplot for the TDRs using the intervention analysis based on the same 10 sets of simulation. It shows that the variation among different sets of simulation is small. However, the median TDR is around 0.1. This shows that it is hard to find the imputed anomalies by using the intervention analysis. Note that the intervention analysis method actually treats the background time series as an ARMA sequence, which is a wrong model in our case. The bottom panel shows the naive one-scale detection method, which is discussed in Section 3.2. The median TDR for this method is larger than the intervention analysis, but smaller than the MRAD method. Note that this naive method also assumes a wrong model, that all the observations are independent. These three plots suggest that our MRAD is (on average) the best among these three methods.

The middle panel and the right panel in Figure 4 shows the distribution of the FDRs and FNRs for these three methods using the same 10 sets of simulation. The left panel is the boxplot of the FDRs, and the right panel shows the boxplot of the FNRs. By comparing these three methods, we find that the median FDR of the MRAD method is the smallest one, and the intervention models has the largest median FDR. All three methods have really small median FNRs. In fact, most of these methods report a large proportion of the data observations as regular observations. This will cause the FNRs to be relatively small. In addition, we notice that the median FNR of the MRAD method is the smallest as well, which suggests that the MRAD method provides a better solution in detecting outliers, when the background time series is a long range dependent trace.

More simulation results are reported in Zhang (2007), in which we imputed the local mean level shift with different combinations of durations, starting points and intensities. It is shown that our method is more effective in identifying the imputed anomalies, even when the intensity is small (e.g., around 1/2 of the standard deviation of the background trace). We will not report them here in order to save space.

### *5.3. Semi-experiments*

In the field of statistical anomaly detection, there is a lack of labeled data sets, and we can uses semi-experiment to evaluate the performance of the MRAD method. The idea is that one first collects some real traces from the Internet, and filter out abnormal parts to form the "normal" (background) trace. One then simulates some well-known network anomalies, such as port scans, rose attacks and TCP SYN flood attacks, and combines these two types of traces together, to form a testbed for anomaly detection methods. Here we briefly show one semi-experiment, in which we find that our method is more efficient when the imputed anomalies have low intensity with long duration. This example has been discussed in Zhang et al. (2008) as well. See more detailed analysis of this example in Zhang (2007).

A three-hour real network trace was collected from the UNC campus Internet. We removed those network flows without either the starting point or the end point within the trace, and treated the remaining flows as the "normal" traffic. Because of this special treatment, the data at the beginning and the end of this trace will have relatively small magnitude than the central part. In order to obtain an approximately stationary time series, we use the central one-hour trace as the background traffic. In this semi-experiment, we show one example of port scan as the injected anomaly.

A port scan generates a series of small packets to learn which computer network services ( associated with a port number), the target computer provides. The port scan can give the attackers information as to where to probe for system weaknesses. Usually a port scan with high frequency significantly increases the number of flows (i.e., a huge local mean level shift). It also increases the number of packets and bytes to a certain degree. If the sizes (in bytes) of the probing packets are really small, the increase in byte counts will typically be dominated by the variability of the real trace itself. In addition if the probing uses a low frequency, the change of packet counts will also be dominated by the variability of the background. In this situation, the port scans are not detectable in the series of packet counts and byte counts. For this example, we used a medium frequency to send out small probing packets, such that the increase in packet count time series is detectable. See Lee et al. (2003) for more about detection methods and characterization of port scan attacks.

Note that because the "normal" trace is generated from a collected real trace, the MRAD method might also flag some anomalies already presented in the background. After combining the anomaly trace with the background trace, we computed the packet, byte, and flow counts per $10ms$ time interval. As discussed earlier, the background trace used here lasts for 1 hour. The port scan simulated for this example lasted 6 minutes, i.e., 10% of the total trace.

Figure 5 (a.k.a. Figure 4 in Zhang et al. (2008)) shows the time series of the packet-count trace in the top panel. It is hard to tell whether there are network anomalies, and to identify the locations of the anomalies. This time series has the estimated Hurst parameter as 0.95, which indicates LRD. We use this estimate (0.95), to perform an MRAD procedure. The bottom panel in Figure 5 is the

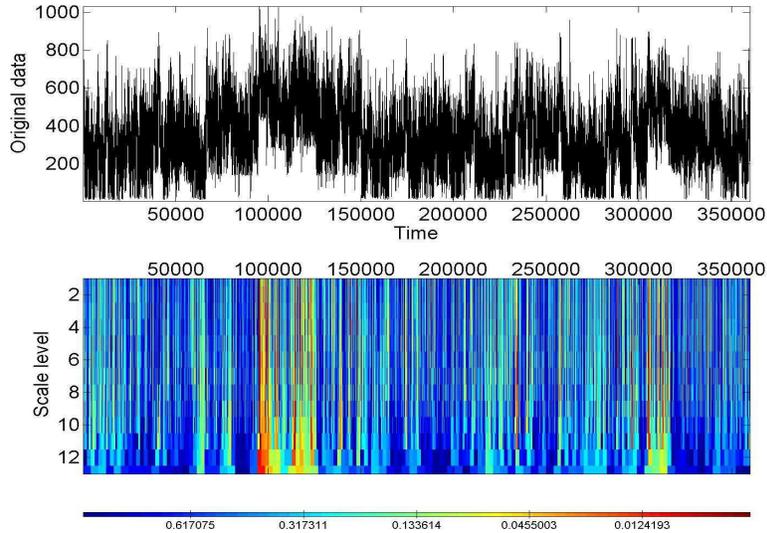

Fɪɢ 5. *The color MRAD outlier map for the packet-count series of the semi-experiment example. It highlights two possible outlying regions: around 100,000, and around 300,000.*

MRAD outlier map based on NOWA. Two hot zones are highlighted in the map: around 100,000 and around 300,000. The zone around 100,000 corresponds to the locations where we injected the network anomaly. We conjecture that those anomalies around 300,000 are some actual anomalies within the original Internet trace, and we are working on identifying these anomalies.

## 6. Discussion

We have shown that the MRAD method is more efficient for identifying outliers in LRD time series, compared with single scale methods, especially when outliers are in the form of a slight mean level shift, with relatively small intensity. We also proved that for a two-scale procedure, the MRAD method uses a more conservative threshold, and has larger power on average than detection methods based on a single scale.

In this paper a naive aggregation method is used to illustrate the advantage of multi-resolution methods and prove theoretical results. In practice other aggregation methods such as kernel methods and wavelet methods can be used to form the multi-resolution time series, which may lead to more powerful detection algorithms. This paper discussed one type of outlier (local mean level shift) which is very common in the DDoS attack, while detection of other types of outliers such as those in Tsay (1988) from a LRD time series is of both theo-retical and practical importance. These are questions we plan to address in our future research.

As discussed in Section 1, using wavelet-based change point detection methods can also help to identify network anomalies. Examining the connections and the differences between our methods and the change point detection methods is potential future work. In addition, the method described in Wang (1999) provided a way to incorporate the multiple comparisons in the time domain at a given scale. It will be helpful to develop a new method based on our method or their method, which uses multiple comparisons in both the scale and the time spaces. This adjustment will improve the performance of the detector.

## Acknowledgement

This work is part of the first author's PhD dissertation under the supervisor of the second and third authors. The authors thank Jeff Terrell for generating the semi-experiments. Thanks also go to Professor Richard Smith, for discussion on extreme value theory. This work was partially supported by NSF fund DMS-0308331 and DMS-0605434.

## Appendix

**Proof for Proposition 1**

Due to the fact that $\{|Y_L(i)| > C_\alpha^M\} \subset \bigcup_L \{|Y_L(i)| > C_\alpha^M\}$, we have

$$P\{|Y_L(i)| > C_\alpha^M\} \leq P\left(\bigcup_L \{|Y_L(i)| > C_\alpha^M\}\right) = P\{|Y_L(i)| > C_\alpha\},$$

which yields that $C_\alpha \leq C_\alpha^M$.

**Proof for Proposition 2**

Let $b$ the aggregation bin size, and $l = b, b^2, \cdots, b^k, \cdots$ are the aggregation scales. Define $\xi_b(k) = Y_{b^k}(1)$. Let $L = b^{k+1}$, Lemma 3.8.3 (page 123) in Zhang (2007) gives

$$\rho_b(k) = \text{Cov}(\xi_1, \xi_{k+1}) = \text{Cov}(X_1(1), X_L(1)) = \frac{1}{b^{kH}}\left[1 + \frac{1}{2}(b^{k2H} - (b^k - 1)^{2H} - 1)\right]$$

$$= Hb^{k(H-1)} + \frac{b^{-kH}}{2} - \frac{H(2H-1)}{2}b^{k(H-2)} + o(b^{k(H-2)}),$$

which is the proposition.

**Proof for Proposition 3**

Assume the correlation coefficient between $Y_L(i)$ and $Y_1(i)$ is $\rho$. When $H_0$ is true, we have the following result, by using the bivariate normal calculation (e.g., Lemma 3.8.1 in Zhang (2007) and its remarks).

$$P_0(\{|Y_L(i)| > C_\alpha^M\} \cup \{|Y_1(i)| > C_\alpha^M\}) = 1 - P_0(\{|Y_L(i)| \leq C_\alpha^M\} \cap \{|Y_1(i)| \leq C_\alpha^M\})$$
$$= 1 - [(\Phi(C_\alpha^M) - \Phi(-C_\alpha^M))^2 + 2(C_\alpha^M)^2 \phi^2(C_\alpha^M)\rho^2 + O(\rho^3)].$$

Let $L \to \infty$, we have $\rho \to 0$, the above converges to $1 - (\Phi(C_\alpha^M) - \Phi(-C_\alpha^M))^2 = 1 - (1 - 2\Phi(-C_\alpha^M))^2 = \alpha$, which leads to $C_\alpha^M = -\Phi^{-1}((1 - \sqrt{1 - \alpha})/2)$.

Let $\rho = \text{Cov}(Y_1(i), Y_L(i))$. Lemma 3.8.3 in Zhang (2007) gives $\rho = HL^{(H-1)} + o(L^{(H-1)})$. In addition, let $C_\alpha^M = C_0 + a\rho^\gamma + o(\rho^\gamma)$. When $\rho \to 0$, we have

$$P_0(\{|Y_L(i)| > C_\alpha^M\} \cup \{|Y_1(i)| > C_\alpha^M\}) = 1 - [(\Phi(C_\alpha^M) - \Phi(-C_\alpha^M))^2 + 2(C_\alpha^M)^2 \phi^2(C_\alpha^M)\rho^2 + O(\rho^3)].$$

Using Taylor expansion as $\rho \to 0$, We have

$$P_0(\{|Y_L(i)| > C_\alpha^M\} \cup \{|Y_1(i)| > C_\alpha^M\}) = 1 - P_0(\{|Y_L(i)| \leq C_\alpha^M\} \cap \{|Y_1(i)| \leq C_\alpha^M\})$$
$$= 1 - [(\Phi(C_\alpha^M) - \Phi(-C_\alpha^M))^2 + 2(C_\alpha^M)^2 \phi^2(C_\alpha^M)\rho^2 + O(\rho^3)]$$
$$= \alpha - 4\phi(C_0)a\sqrt{1 - \alpha}\rho^\gamma - 2\phi^2(C_0)C_0^2\rho^2 + O(\rho^2). \tag{3}$$

(See detail of (3) in page 126 of Zhang (2007)).

From the above, we have $\gamma = 2$, and $a = -\phi(C_0)C_0^2/(2\sqrt{1 - \alpha})$.

In summary, we have

$$C_\alpha^M = C_0 - \frac{\phi(C_0)C_0^2}{2\sqrt{1 - \alpha}}\rho^2 + o(\rho^2).$$

Substituting with $\rho = HL^{H-1} + o(L^{H-1})$, we have the proposition.

**Power for a two-scale MRAD procedure**

**Lemma 1.** *Let* $\beta_{(1,L)} = P_1(\max_{l=1,L}\{|Y_l(i)| > C_\alpha^M\})$, *we have*

$$\beta_{(1,L)} = 1 - \sqrt{1 - \alpha}[\Phi(C_0 - \delta) - \Phi(-C_0 - \delta)] - 2Hk\delta C_0\phi(C_0)[\phi(C_0 - \delta) - \phi(-C_0 - \delta)]L^{-1} + o(L^{-1}).$$

*when* $L \to 0$. *Here* $C_0 = \Phi^{-1}((1 + \sqrt{1 - \alpha})/2)$.

The proof of this lemma uses the bivariate normal calculation and the above proposition. See page 128 of Zhang (2007) for details.

**Proof for Theorem 1**

Let $\rho = \text{Cov}(Y_1(i), Y_L(i))$. From the equation (1), we know when the alternative hypothesis is true, we have the marginal distribution of $Y_1(i)$, and $Y_L(i)$ are $Y_1(i) \sim N(\delta, 1)$, $Y_L(i) \sim N(\mu_L, 1)$, where $\mu_L = K\delta/L^H$, and $K = 1, 2, \cdots, L$.

From the Lemma 3.8.2 (page 122) in Zhang (2007), the power at scale 1 and $L$ is given by

$$\beta_1 = 1 - [\Phi(C_\alpha - \delta) - \Phi(-C_\alpha - \delta)], \quad \beta_L = 1 - [\Phi(C_\alpha - \mu_L) - \Phi(-C_\alpha - \mu_L)].$$

When $K$ and $\delta$ are fixed, let $L \to \infty$, we have $\rho \to 0$, $\mu_L \to 0$, and

$$\beta_L = \alpha + O(L^{-2H}), \quad \beta_{(1,L)} = 1 - [\Phi(C_0 - \delta) - \Phi(-C_0 - \delta)]\sqrt{1-\alpha} + O(L^{-1}),$$

which can be directly derived from Lemma 3.8.2 in Zhang (2007) and Lemma 1.

Define the power difference function $f(\alpha, \delta)$ as

$$f(\alpha, \delta) = \left[\Phi\left(\Phi^{-1}\left(1 - \frac{\alpha}{2}\right) - \delta\right) - \Phi\left(-\Phi^{-1}\left(1 - \frac{\alpha}{2}\right) - \delta\right)\right] + [1 - \alpha]$$
$$- 2\sqrt{1-\alpha}\left[\Phi\left(\Phi^{-1}\left(\frac{1+\sqrt{1-\alpha}}{2}\right) - \delta\right) - \Phi\left(-\Phi^{-1}\left(\frac{1+\sqrt{1-\alpha}}{2}\right) - \delta\right)\right].$$

We have $2\beta_{(1,L)} - (\beta_1 + \beta_L) = f(\alpha, \delta) + O(L^{-1})$, when $\alpha \to 0$. Thus we only need to show $f(\alpha, \delta) > 0$ as $\alpha \to 0$.

Note that $\alpha = 0$, $f(\alpha, \delta) = 2 - 2 = 0$.

$$\frac{\partial f(\alpha, \delta)}{\partial \alpha} = -\frac{1}{2}\exp\left\{-\frac{\delta^2}{2}\right\}[\exp\{\delta C_\alpha\} + \exp\{-\delta C_\alpha\}] - 1 + \frac{1}{\sqrt{1-\alpha}}[\Phi(C_0 - \delta) - \Phi(-C_0 - \delta)]$$
$$+ \frac{1}{2}\exp\{-\frac{\delta^2}{2}\}[\exp\{\delta C_0\} + \exp\{-\delta C_0\}].$$

By the *mean value theorem* in Calculus (e.g., page 43 in Beals (1973)),

$$C_0 - C_\alpha = \Phi^{-1}\left(\frac{1+\sqrt{1-\alpha}}{2}\right) - \Phi^{-1}\left(1 - \frac{\alpha}{2}\right)$$
$$= \frac{1}{\phi(\Phi^{-1}(\xi))}\left[\frac{\sqrt{1-\alpha} - (1-\alpha)}{2}\right], \quad (4)$$

where $\xi \in (1 - \alpha/2, (1 + \sqrt{1-\alpha})/2)$.

Equation (4) implies

$$\frac{1}{\phi(C_\alpha)}\left[\frac{\sqrt{1-\alpha} - (1-\alpha)}{2}\right] \leq C_0 - C_\alpha \leq \frac{1}{\phi(C_0)}\left[\frac{\sqrt{1-\alpha} - (1-\alpha)}{2}\right]. \quad (5)$$

Since

$$\lim_{\alpha \to 0+} \frac{1}{\phi(C_0)}\left[\frac{\sqrt{1-\alpha} - (1-\alpha)}{2}\right] = \lim_{\alpha \to 0+} \frac{\frac{1}{2} \times [\frac{1}{2}(1-\alpha)^{-1/2} \times (-1) + 1]}{\Phi^{-1}(\frac{1+\sqrt{1-\alpha}}{2}) \times \frac{1}{4} \times (1-\alpha)^{-1/2}} = 0,$$
$$\lim_{\alpha \to 0+} \frac{1}{\phi(C_\alpha)}\left[\frac{\sqrt{1-\alpha} - (1-\alpha)}{2}\right] = \lim_{\alpha \to 0+} \frac{\frac{1}{2} \times [\frac{1}{2}(1-\alpha)^{-1/2} \times (-1) + 1]}{\Phi^{-1}(1 - \frac{\alpha}{2}) \times \frac{1}{2}} = 0,$$

we have $\lim_{\alpha \to 0+}(C_0 - C_\alpha) = 0$.

Thus, when $\alpha \to 0$, $\exp\{\delta(C_0 - C_\alpha)\} - 1 = \delta(C_0 - C_\alpha) + o(\delta(C_0 - C_\alpha))$, which leads to

$$\exp\{\delta C_0\} - \exp\{\delta C_\alpha\} = \exp\{\delta C_\alpha\}\exp\{\delta(C_0 - C_\alpha) - 1\}$$
$$= \exp\{\delta C_\alpha\}[\delta(C_0 - C_\alpha) + o((C_0 - C_\alpha))].$$

By equation (5),

$$\lim_{\alpha \to 0+}\exp\{\delta C_\alpha\}(C_0 - C_\alpha) \geq \lim_{\alpha \to 0+}\frac{(\sqrt{1-\alpha} - (1-\alpha))/2}{\exp\{-\delta C_\alpha\}\phi(C_\alpha)}$$
$$= \lim_{\alpha \to 0+}\frac{\frac{1}{2} \times [1 - \frac{1}{2}(1-\alpha)^{-1/2}]}{\exp\{-\delta C_\alpha\}[\frac{\delta + C_\alpha}{2}]} = +\infty.$$

This yields that $\exp\{\delta C_0\} - \exp\{\delta C_\alpha\} \to +\infty$, when $\alpha \to 0+$, i.e.

$$\frac{\partial f(\alpha, \delta)}{\partial \alpha}|_{\alpha \to 0+} > 0.$$

From the above, we know that for any $\delta > 0$, there exists $\alpha_\delta$, such that for any $\alpha \in (0, \alpha_\delta)$, we have $f(\alpha, \delta) > 0$, hence the theorem holds.